\documentclass[12pt]{amsart}

\date{\today}

\usepackage[margin=1.15in]{geometry}

\usepackage{amsmath,amscd,amssymb,amsfonts,latexsym,wasysym, mathrsfs, enumitem, mathtools,hhline,color}
\usepackage[all, cmtip]{xy}

\usepackage{url}

\definecolor{hot}{RGB}{65,105,225}

\usepackage[pagebackref=true,colorlinks=true, linkcolor=hot ,  citecolor=hot, urlcolor=hot]{hyperref}%make all the references and links clickable

\usepackage{graphicx,enumerate}

\newcommand{\C}{\mathbb{C}}

\theoremstyle{plain}
\newtheorem{theorem}{Theorem}[section]
\newtheorem{prop}[theorem]{Proposition}

\newtheorem{lm}[theorem]{Lemma}

\newtheorem{cor}[theorem]{Corollary}

\newtheorem{thrm}[theorem]{Theorem}

\theoremstyle{definition}

\newtheorem{defn}[theorem]{Definition}

\newtheorem{rmk}[theorem]{Remark}

\newtheorem{ex}[theorem]{Example}
\newtheorem*{ex*}{Example}

\def\be{\begin{equation}}
\def\ee{\end{equation}}

\def\bt{\begin{thrm}}
\def\et{\end{thrm}}

\def\bc{\begin{cor}}
\def\ec{\end{cor}}

\def\br{\begin{rmk}}
\def\er{\end{rmk}}

\def\bp{\begin{prop}}
\def\ep{\end{prop}}

\def\bl{\begin{lm}}
\def\el{\end{lm}}

\def\bex{\begin{ex}}
\def\eex{\end{ex}}

\def\bd{\begin{defn}}
\def\ed{\end{defn}}

\newcommand\sV{{\mathcal V}}
\newcommand\sO{\mathcal{O}}

\newcommand{\fJ}{\mathbf{J}}
\newcommand{\fL}{\mathbf{L}}

\newcommand\first{{\mathbf{1}}}

\DeclareMathOperator{\ord}{ord}

\DeclareMathOperator{\cha}{Char}

\def\bC{\mathbb{C}}

\def\bQ{\mathbb{Q}}

\def\bZ{\mathbb{Z}}

\title{The Deletion Restriction method for plane curve germ complement}

\author{Yongqiang Liu}
\address{The Institute of Geometry and Physics, University of Science and Technology of China, 96 Jinzhai Road, Hefei 230026 P.R. China}
\email{liuyq@ustc.edu.cn}



\date{\today}

\begin{document}

\maketitle

\begin{abstract}  We study the cohomology jump loci of rank one local systems for the plane curve germ complement by the deletion restriction method. In particular, the deletion restriction type results in this case only depend on the linking numbers between the irreducible components of the plane curve germ.
 \end{abstract}
 
%%%%%%%%%%%%55 
\section{Introduction}
Let $f: (\bC^2,0) \to (\bC, 0 ) $ be a reduced plane curve singularity at the origin of $\bC^2$. Let  $f=\prod_{i=1}^r f_i$ be a decomposition of the germ $f$ into irreducible factors $f_i$ in the factorial ring $\sO_2$ of analytic function germs. The assumption $f$ being reduced means that all the factors $f_i$ above are distinct. 
   
  In this paper, we always assume that $r\geq 2$.  Let $B$ denote the small enough open ball centred at the origin of $\bC^2$.   Set \begin{center}
    $C_i=\{f_i=0\} \cap B$, $U=B- \bigcup_{i=1}^r C_i$ and $U'= B-\bigcup_{i=2}^r C_i$. 
  \end{center}
Let $j:U\to U'$ be the corresponding inclusion map, which induces a natural map on the fundamental group level: $$ \pi_1(U) \to \pi_1(U').$$

  What can we say about this natural map?  In this paper, we answer the question on the cohomology jump loci level. Since the homotopy type of the local complement are determined by Puiseux pairs and linking numbers (see \cite[Theorem 2.2.11]{D2}), hence so is this natural map. It turns out that our answer on the cohomology jump loci level  only depends on the linking numbers between $C_1$ and $C_i$ for all $1<i \leq r$.

\medskip

First let us recall the definition of cohomology jump loci. %hence the moduli space of rank one $\bC$-local system is $M_B(U)\cong (\bC^*)^r$.
Let $\cha(U)$ be the moduli space of rank one $\bC$-local systems on $U$. %Note that $H_1(U,\bZ)=\bZ^r$. 
Then  $\cha(U)\cong \mathrm{ Hom} (H_{1}(U, \bZ),\C^{\ast})$ as algebraic group. Notice that $H_1(U, \bZ)$ has a natural basis $[\gamma_i]_{1\leq i \leq r}$, where each $\gamma_i$ is the meridian around the irreducible component $C_i$. This basis induces a natural coordinate system $$\cha(U)\cong (\C^{\ast})^{r}=\{t=(t_1, \ldots, t_r)\mid t_i\in \bC^* \text{ for any } 1\leq i\leq r\}.$$
%For any $t=(t_{1},\cdots, t_{r}) \in (\C^{\ast})^{r}$, let $L_t$ be the corresponding rank one local system on $U$, which sends each meridian associated to $H_{i}$ in $\pi_{1}(U)$ to the non-zero complex number $t_{i}$.
 For any $t\in \cha(U)$, let $L_t$ denote the corresponding rank one $\bC$-local system on $U$. %$t$ is sometime omitted, when the context is clear.
\bd  \label{jump loci} The {\it cohomology jump loci of $U$} are defined as:
$$ \sV^{i}_k(U): =  \{t \in \cha(U) \mid   \dim H^i(U,  L_t) \geq k\}.
$$ 
% where $L_t$ is as before the rank-one $\bC$-local system on $U$ associated to $t\in M_B(U)$.
\ed
It is  well known that $\sV^i_k(U)$ are closed sub-varieties of $\cha(U)$ and  homotopy invariants of $U$. In particular, $\sV^1_k(U)$ only depends on the fundamental group $\pi_1(U)$.

The open inclusion map $j:U\to U'$ induces an embedding on the moduli space
$$ j^\#: \cha(U')    \hookrightarrow \cha(U) ,$$ which sends $
(t_2,\cdots,t_r) \in \cha(U') $ to $ (1,t_2,\cdots,t_r)\in \cha(U).$

 \bt \label{main} With the above assumptions and notations, we have that for any $k\geq 1$  
 $$j^\#( \sV^1_k(U'))=\{t_1=1\} \cap \bigg(\sV^1_{k+1}(U) \cup \Big( \sV^1_{k} (U) \cap \{\prod_{i=2}^r t_i^{l_{1i}}\neq 1 \}\Big) \bigg),$$
 where $l_{1i}>0$ is the linking number between $C_1$ and $C_i$. 
 \et 
 
  \br 
Note that $H^0(U, L_t)\neq 0$ if and only if $L_t$ is the trivial local system, hence
$\sV^0_1(U)=\{\first\}$. Here $\first$ denotes the trivial rank one $\bC$-local system (i.e., the constant sheaf) on $U$. 
Since $\chi(U)=0$ and $U$ is homotopy equivalent to a 2-dimensional CW complex, we have that $\dim H^1(U, L_t) =\dim H^2(U, L_t)$ for any $L_t\neq \first$. So for any positive integer $k$, we have that $$\sV^1_k (U) \setminus \{\first\}= \sV^2_k(U) \setminus \{\first\}.$$ The same claim also holds for $U'$. Hence one can get a similar formula for $\sV^2_k$ as in Theorem \ref{main}. In particular, this shows that one can read the jump loci of $U'$ from the one for $U$, once the linking numbers $l_{1i}$ (for $2\leq i \leq r$) are all known.
\er

\br \label{equality} Let $L_t$ be a rank one $\C$-local system on $U$ with $t_1=1$. 
 Let $L'_t$ be the unique rank one $\bC$-local system on $U'$ such that $j^{-1}L'_t=L_t$.
Theorem \ref{main} holds, if we can show that   \be \label{key}
\dim H^1(U', L'_t)=  \left\{ \begin{array}{ll}
\dim H^1(U, L_t),  &  \text{ if } \prod_{i=2}^r t_i^{l_{1i}}\neq 1, \\
\dim H^1(U, L_t) -1, &  \text{ if } \prod_{i=2}^r t_i^{l_{1i}}=1.\\
\end{array}\right. 
\ee
The proof is based on the deletion-restriction method, which is a fundamental tool in the theory of hyperplane arrangements, see \cite[Theorem 3.5]{D3}.
\er
 \br Let $L_t$ be a non-trivial rank one $\C$-local system on $U$ with $t_i=1$ for some $i$.  Without loss of generalities, we assume that $t_1=\cdots=t_s=1$ and $t_i \neq 1$ for any $s+1 \leq i \leq r$. Set $V= B- \bigcup_{i=s+1}^r C_i$ and let $J_t$ be the unique rank one local system on $V$ such that $J_t\vert_{U}=L_t$.  Using formula (\ref{key}) repetitively, one can always reduce the computation of $H^1(U, L_t)$ to $H^1(V, J_t)$, once the linking numbers are known.  For more results in this direction, see Theorem \ref{split}.
 \er

 \textbf{Acknowledgements.} The author is grateful to Javier Fernandez de Bobadilla, Nero Budur, Jos\'{e} Ignacio Cogolludo Agust\'{i}n  and Alexandru Dimca for useful discussions.  %The author also would like to thank Yuejiao Li for translating some Russian references.

 %For the torsion local system, we can work in the category of mixed Hodge modules. A claim due to P. Brosnan, H. Fang, Z. Nie and G. Pearlstein \cite[Lemma 2.18]{BFNP} shows that the intermediate extension functor is exact upon certain weight conditions. We finish the proof of Theorem  by verifying these weight conditions.

%%%%%%%%%%%%%%%%%%%%%%%%%%%%%%

\section{The deletion restriction method}

In this section, we always assume that the rank one local system $L_t$ has $t_1=1$.
  Note that $C_1^*=U'-U=C_1-\{0\}$ is homotopy equivalent to the real circle $S^1$. Define a rank one $\bC$-local system $L''_t(1)$  on $C_1^*$, whose representation sends the generator  of the fundamental group of $C_1^*$ to $\prod_{i=2}^r t_i^{l_{1i}}$.  
 
 \bl \label{dr} For any $t\in \cha(U)$ with $t_1=1$, we have that $(R^1 j_* L_t)\vert _{C_1^*}= L''_t(1).  $ 
Moreover, we have a short exact sequence of perverse sheaves:
\be 
 \label{ses1}
 0 \to L'_t[2] \to Rj_* L_t[2] \to L''_t(1)[1] \to 0 
\ee Here $L''_t(1)[1]$ is considered as a perverse sheaf on $U'$ with support on $C_1^*$.
 \el
\begin{proof} Note that $j: U \to U'$ is a Stein morphism, hence $j_!L_t[2]$ and $Rj_*L_t[2]$ are both perverse sheaves.  Let $j_{!*}$ denote the intermediary extension functor, e.g., see \cite[Definition 5.2.6]{D2}. 
Consider the short exact sequence of perverse sheaves:
$$0 \to P \to j_! L_t[2] \to j_{!*} (L_t[2]) \to 0  ,$$ 
 where $P$ is the kernel. Note that $j_{!*} (L_t[2])=L'_t[2]$. Meanwhile, $P$ is supported on $C_1^*$ and $$ P\vert_{C_1^*}= L'_t[2]\vert_{C_1^*} [-1] = L'_t\vert_{C_1^*} [1].$$
 
 Note that the inclusion from $C_1^* $ to $U'$ induces a map $\rho: H_1(C_1^*,\bZ) \to H_1(U',\bZ)$. Let $\gamma_i$ denote the meridian associated to $C_i$ in $U'$. Let $\gamma$ denote the generator of $\pi_1(C_1^*)$.   Then $\rho ([\gamma])=\sum_{i=2}^r l_{1i} \cdot [\gamma_i]$. In fact, $C_i^*=C_i -\{0\}$ is homotopy equivalent to a knot, denoted by $K_i$ for any $1\leq i\leq r$.   Consider the Seifert surface $\Sigma_1$ of the knot $K_1$. Then $\Sigma_1$ intersects with the knot $K_i$ ($i\geq 2$) at $l_{1i}$ points  due to the definition of linking numbers, e.g. see \cite[Definition 2.1.7, Proposition 2.2.12]{D1}. Then the map $\rho$ can be viewed as the composition of the following maps:
\be \label{comp}  H_1(K_1,\bZ) \to H_1(\Sigma_1 \setminus \cup_{i=2}^r K_i, \bZ) \to H_1(U',\bZ) .\ee
 For any point $p\in \Sigma_1\cap K_i$ ($i\geq 2$), let $\gamma_p$ denote the meridian around $p$ in $\Sigma_1$. Then the last map in (\ref{comp}) sends $[\gamma_p]$  to $[\gamma_i]$ for any  $p\in \Sigma_1\cap K_i$ and the first map in (\ref{comp}) sends $[\gamma]$ to $ \sum_{p} [\gamma_p]$, where the sum runs over all points $p \in \Sigma_1\cap (\cup_{i=2}^r K_i)$. %Then $\rho ([\gamma])=\sum_{i=2}^r l_{1i} \cdot [\gamma_i]$ follows.
 
 So  $ P\vert_{C_1^*} = L'_t\vert_{C_1^*} [1]=L''_t(1)[1].$ We rewrite the short exact sequence as follows:
 $$
 0\to L''_t(1)[1] \to  j_! L_t[2]  \to L'_t[2] \to 0
 $$ Consider  Verdier dual of this short exact sequence:
$$
 0 \to (L'_t)^\vee[2] \to Rj_* L^\vee_t[2] \to (L''_t(1))^ \vee[1] \to 0. 
 $$ Here $L_t ^\vee$ denote the corresponding dual local system.
Hence $(R^1 j_* L^\vee_t)\vert_{C_1^*}= (L''_t(1))^\vee$.
  Note that $L^\vee_t= L_{t^{-1}}$, where $t^{-1}=(t_1^{-1},\cdots,t_r^{-1})$. Then the claims follow by replacing $t$ by $t^{-1}$.    
\end{proof}
  %Set $U'= X-\bigcup_{i=2}^r V_i$, and let $i:U\to U'$ be the inclusion map.  Let $L'$ be the rank one $\bC$-local system on $U'$ such that $i_{!*}(L[2])=L'[2]$. Note that $V_1^*=U'-U$ is homotopy equivalent to the real circle. Define the rank one local system $L''$ on $V_1^*$, whose representation sends the only generator  of the fundamental group to $\prod_{i=2}^r t_i$.  

\br
The short exact sequence (\ref{ses1}) is in the same spirit as the deletion-restriction method used in hyperplane arrangements, e.g. see \cite[page 221]{D2}. It induces a hypercohomology long exact sequence:
\be  \label{les}
\begin{split}
0 \to H^1(U',L'_t) \to H^1(U,L_t)  \to H^0(C_1^*, L_t''(1) ) \overset{\delta}{\to}  H^2(U',L'_t) \\
\to H^2(U,L_t) \to H^1(C_1^*, L_t'' (1)) \to 0 
\end{split}
\ee
In the proof of Theorem \ref{main}, we will show that the boundary map $\delta$ in (\ref{les}) is always trivial.

A long exact sequence similar to (\ref{les}) in  the hyperplane arrangement case was firstly obtained by D. C. Cohen in \cite{Co}. Note  that in this case there are examples where the boundary map is non-trivial for
some rank one $\C$-local system, see \cite[Example 4.4]{BL}.
\er

\section{Proof of Theorem \ref{main}}
We first recall the following structure theorem due to N. Budur and B. Wang \cite{BW2}.
\bt[Structure Theorem] Let $U$ be the complement in a small ball of a germ of a complex analytic set. Then 
$\sV^i_k(U)$ is a finite union of torsion translated subtori of $\cha(U)$ for any non-negative integers $i$ and $k$.
\et

Fix an integer $s$ such that $1\leq s <r$. Let $L_t$ be a rank one local system on $U$ with  $t_1=\cdots=t_s=1$ and $t_i \neq 1$ for any $s+1 \leq i \leq r$. Set $V= B- \bigcup_{i=s+1}^r C_i$ and let $u: U\to V$ denote the inclusion map from $U$ to $V$.
  Let $J_t$ be the unique rank one local system on $V$ such that $u^{-1}J_t =L_t $.
By Lemma \ref{dr}, we have the following short exact sequences of perverse sheaves:
\be \label{ses2} 
0 \to J_t[2] \to Ru_* L_t[2] \to \oplus_{i=1}^{ s} L''_t(i) [1] \to 0 ,
\ee
where $L''_t(i)$ is the rank one local system supported on $C_i^*$ as defined at the beginning of the previous section.  
Then it induces the following long exact sequence:
\be \label{split2}
\begin{split}
0 \to H^1(V,J_t) \to H^1(U,L_t)  \to \oplus_{i=1}^s H^0(C_i^*, L_t''(i) ) \overset{\partial}{\to}  H^2(V,J_t) \\ \to H^2(U,L_t) \to \oplus_{i=1}^s H^1(C_i^*, L_t''(i) ) \to 0 
\end{split}
\ee

\bt \label{split} With the above assumptions and notations, the boundary map $\partial$ in the long exact sequence (\ref{split2}) is always trivial.
\et 
\begin{proof}
%Proving Theorem \ref{main} is equivalent to show that the boundary map $\partial$ in the above long exact sequence  is always trivial.  
The proof is divided into 3 steps. 

\medskip

\item[Step 1:] We first reduce the proof to the torsion local system case by the structure theorem.

\medskip 

By an {\it absolute $\bQ$-constructible subset} (see \cite[Theorem 1.3.1]{BW18}) of the space of rank one local systems $\cha(U)$, we will mean a Zariski constructible subset defined over $\bQ$, which is obtained from finitely many torsion-translated complex affine algebraic subtori of $\cha(U)$ via a finite sequence of taking union, intersection, and complement. 

Let us fix some notations: $t_{\leq s} \equiv 1$ stands for that $t_i=1$ for any $1\leq i \leq s$ and $t_{>s} \neq 1$ stands for that  $t_i\neq 1$ for any  $s+1\leq i \leq r$.  
Consider the set $$\Theta:=\{t\in    \cha(U)\mid t_{\leq s} \equiv 1, t_{>s}\neq 1 \text{  and the boundary map } \partial=0  \text{  for } L_t \}.$$ 

Note that  $$\dim H^0(C_1^*, L''_t(1))= \left\{ \begin{array}{ll}
0,  &  \text{ if } \prod_{i=2}^r t_i^{l_{1i}}\neq 1, \\
1, &  \text{ if } \prod_{i=2}^r t_i^{l_{1i}}=1.\\
\end{array}\right.  $$ Similar result also holds for $L''_t(i)$ with $1\leq i \leq s$. 
Together with the structure theorem for the jump loci of $U$ and $V$, it implies that  $\Theta$ is an absolute $\bQ$-constructible subset of $\cha(U)$.   
This reduces the proof to the torsion local system case. Indeed, if the boundary map $\partial=0$ holds for any torsion local system $L_t$ with $t_{\leq s}\equiv 1  $ and $t_{>s}\neq 1$, then the smallest  absolute $\bQ$-constructible set containing all these torsion local systems is exactly $ \{t\in \cha(U)\mid  t_{\leq s} \equiv 1, t_{>s}\neq 1 \}$.
 So it implies $\partial=0$ for any $L_t$ with $ t_{\leq s} \equiv 1$ and $ t_{>s}\neq 1$.

From now on, we always assume that $L_t$ with $ t_{\leq s} \equiv 1$ and $ t_{>s}\neq 1$ is a torsion local system, that is, $t_i\neq 1$  are roots of unity for any $s+1\leq i\leq r$.  %It is easy to see that (\ref{key}) holds for the constant sheaf case, so we assume that $t_i\neq 1$ for some $i\geq 2$.

\medskip

\item[Step 2:] We prove a $\bQ$-version of the short exact sequence (\ref{ses2}), which is indeed a short exact sequence of mixed Hodge modules. Then we extend  it to $B$ by the intermediary extension functor.

\medskip 

%The commutative diagram of inclusions 
%\begin{center}
%$\xymatrix{ 
%U' \ar[r] \ar[d]^{j}  & \U' \ar[d]^{\fj}   \\
%U \ar[r] & \U  
%}$
%\end{center}
%induces 
%\begin{center}
%$\xymatrix{ 
%M_B(U')    & M_B(\U') \ar[l]_{\simeq}   \\
%M_B(U) \ar[u] & M_B(\U)  \ar[u] \ar[l]_{\simeq}
%}$
%\end{center}

Note that $t_i=1  \text{  for } 1\leq i \leq s$ and $t_i$  are roots of unity for $i\geq s+1$.  
 Let $O_t$ be the orbit of $t$ in $(\bC^*)^r$ under the diagonal action of the Galois group $Aut(\bC/\bQ)$. Then for every $\alpha \in O_t$, $\alpha_i$ is a primitive root of unity of same order as $t_i$ for every $i$. Hence $O_t$ is finite. The direct sum 
$$
\fL_\bC=\oplus_{\alpha \in O_t} L_\alpha
$$
is a higher-rank local system on $U$, and it is an $Aut(\bC/\bQ)$-invariant $\bC$-point of the moduli space of local systems on $U$. Hence there exists an irreducible  $\bQ$-local system $\fL$ such that $\fL_\bC=\fL\otimes_\bQ\bC$. %In particular, $\fL$ is a simple $\bQ$-local system. 
Moreover, $\fL[2]$ underlies a Hodge module of pure weight 2. 
To see this, one can consider the finite cover $U^t$ of $U$ associated to the following composition of surjective maps $$ \pi_1(U)\to H_1(U,\bZ)\cong \bZ^r \to \oplus_{i=1}^r \bZ/(\ord  t_i) .$$ Here the first map is the Hurewicz morphism,
and the second map sends every factor $\bZ$ to the finite group $\bZ/(\ord t_i)$, where $\ord t_i$ is the order of $t_i$. The finite cover $U^t$ can be taken as a smooth analytic variety such that the covering map $\pi: U^t \to U$ becomes a proper analytic map. Now $\bQ_{U^t} [2]$ underlies a polarized pure  Hodge module of weight 2, since $U^t$ is a complex analytic manifold \cite[Example 15.3]{Sch}. Then so does $R\pi_* \bQ_{U^t} [2]$, since $\pi$ is a finite map. Note that  $\fL[2]$ is a direct summand factor of $R\pi_* \bQ_{U^t} [2]$, so   $\fL[2]$ underlies a polarized  pure Hodge module of  weight 2. In the rest of the proof, we use that all considered functors lift to the category of polarized  mixed Hodge modules, see \cite{Sai}.

 Note that for all $\alpha \in O_t$, $\alpha_i=1$ when $1\leq i \leq s$. Hence we can define a simple $\bQ$-local system $\fJ$ on $V$ in a similar way such that $$\fJ\otimes_{\bQ} \bC= \oplus_{\alpha \in O_t} J_\alpha .$$
  In particular,  $\fJ[2]= u_{!*} \fL[2]$ underlies a polarized Hodge module of weight 2. %Here $\fj$ is the open inclusion from $\U$ to $\U'.$
%Consider the short exact sequence of mixed Hodge modules
%\be \label{algebraic}
%0\to L'[2] \to Rj _* L[2] \to \fQ\to 0 
%\ee

% In the rest of the proof, we always work over the category of mixed Hodge modules with $\bQ$-coefficients.

%Consider the following sequence of maps: 
%$$ U \overset{i }{\hookrightarrow} U' \overset{j' }{\hookrightarrow}    U' \cup \{0\}.$$

 %Note that \cite[Lemma 2.18]{BFNP} is proved for 
 
 %It is shown in  \cite[page 313]{Sai}) that the category of mixed Hodge modules in the algebraic case is equivalent to the category of polarized mixed Hodge modules which are extendable to its closure in the analytic case. In particular, the latter category is preserved by taking restrictions over a open analytic subset.  

%Taking restrictions of (\ref{algebraic}) over the open analytic subset $U$,  we have the following 
Consider the short exact sequence of mixed hodge modules: 
\be \label{ses}
0 \to  \fJ[2] \to   Ru_{\ast}  \fL[2] \to \oplus_{i=1}^s \fL''(i)[1]\mathbf{ (-1)} \to 0    
\ee
where %\begin{center}
%$L\otimes_{\bQ} \bC= \oplus_{\alpha \in O_t} L_\alpha,$
% $L'\otimes_{\bQ} \bC= \oplus_{\alpha \in O_t} L'_\alpha$ and
$\fL''(i)\otimes_{\bQ} \bC= \oplus_{\alpha \in O_t} L''_\alpha(i).$
%\end{center}
The stalk of the local system $\fL''(i)$ comes from the first degree cohomology of a punched disk, and this is where does the Tate twist $\mathbf{ (-1)}$ in (\ref{ses}) come from. Note that the Tate twist $\mathbf{ (-1)}$ changes the weight from 1 to 3, hence $\fL''(i)[1]\mathbf{ (-1)}$ has pure weight 3.

 Consider the inclusion map $v: V \to B$ from $V$ to the small open ball $B$.  With the above weights, \cite[Lemma 2.18]{BFNP} shows that the functor $v_{!*}$ is exact to the short exact sequence (\ref{ses}), hence we get the following short exact sequence of perverse sheaves:
\be \label{newses} 
0 \to  v_{!*}(\fJ[2]) \to  v_{!*} ( Ru_{\ast} \fL[2]) \to v_{!*} \Big(\oplus_{i=1}^s\fL''(i)[1]\mathbf{ (-1)}\Big) \to 0     
\ee
 (\cite[Lemma 2.18]{BFNP} is proved for the category of polarized mixed Hodge modules which are extendable to its closure in the analytic case. According to \cite[Theorem 15.1]{Sch}, the short exact sequence (\ref{newses}) are indeed extendable to the small open ball $B$.)

\medskip

\item[Step 3:] We conclude the proof by showing that the short exact sequence of perverse sheaves (\ref{newses}) implies that the boundary map $\partial$ is trivial.

\medskip

Set $v': V\to B\setminus \{0\}$ and $v'': B\setminus \{0\} \to B$. By the assumption $t_{>s}\neq 1$, we get that $$Rv'_* \fJ [2] =Rv'_! \fJ [2]=v'_{!*} (\fJ[2])$$ and $$ Rv'_* Ru_*\fL [2]=Rv'_!Ru_* \fL [2]=v'_{!*}Ru_*(\fL[2]).$$ 
Recall Deligne's construction for the intermediary extension functor, e.g. see \cite[Proposition 5.2.10]{D2}.  Then the stalk cohomology of $v_{!*}(\fJ[2])= v''_{!*} v'_{!*} (\fJ[2])=v''_{!*}Rv'_* (\fJ [2]) $ at the origin can be computed as follows:  $$H^j(v_{!*}(\fJ[2]))_0 \cong \left\{ \begin{array}{ll}
 H^j(Rv'_* (\fJ [2]))_0 = H^{j+2}(V, \fJ), &  \text{ if } j\leq -1\\
0, & \text{ else}.\\
\end{array}\right. $$
Similarly, we have that $$H^j(v_{!*} Ru_*(\fL[2]))_0 \cong \left\{ \begin{array}{ll}
 H^j(Rv'_* Ru_*(\fL [2]))_0 = H^{j+2}(U, \fL), & \text{ if } j\leq -1\\
0, & \text{ else}.\\
\end{array}\right. $$
Easy computation shows that  $$H^j(v_{!*} \fL''(i)[1])_0 \cong \left\{ \begin{array}{ll}
   H^{j+1}(C_i^*, \fL''(i)), & \text{ if } j\leq -1\\
0, & \text{ else}.\\
\end{array}\right. $$
%It turns out that this claim implies $\delta=0$.

(\ref{newses}) can be viewed as a distinguished triangle. By taking the stalk of (\ref{newses}) at the origin, we get a short exact sequence of $\bQ$-vector spaces: 
$$
0\to H^1(V, \fJ) \to H^1(U, \fL) \to  \oplus_{i=1}^s H^0(C_i^*, \fL''(i)) \to 0
 $$
 It induces the following short exact sequence:
 $$
0\to H^1(V, J_t) \to H^1(U, L_t) \to  \oplus_{i=1}^s H^0(C_i^*, L''_t(i)) \to 0
 $$
 This shows that the boundary map $\partial$ is always trivial.
%Note that   
%$$
%\dim_{\bQ} H^0(C_1^*, \fL'')=  \left\{ \begin{array}{ll}
%0,  &  \prod_{i=2}^r t_i^{l_{1i}}\neq 1, \\
%\# \{O_t\}, &  \prod_{i=2}^r t_i^{l_{1i}}=1.\\
%\end{array}\right. 
%$$ Then the formula (\ref{key}) holds for $\fL$, hence also for $L_t$.
\end{proof}

Now we are ready to prove Theorem \ref{main}.
\begin{proof}[{\bf Proof of Theorem \ref{main}}]
Note that for the local system $L''_t(1)$, $$\dim H^0(C_1^*, L''_t(1))= \left\{ \begin{array}{ll}
0,  &  \text{ if } \prod_{i=2}^r t_i^{l_{1i}}\neq 1, \\
1, &  \text{ if } \prod_{i=2}^r t_i^{l_{1i}}=1.\\
\end{array}\right.  $$
Hence according to Remark \ref{equality}, proving Theorem \ref{main} is equivalent to show that the boundary map $\delta$ in the long exact sequence (\ref{les}) is always trivial. 

If $L_t$ is the constant sheaf, then the claim follows as $b_1(U)=r$ and $b_1(U)=r-1$.

Assume that $L_t$ is not the constant sheaf. Set $s:= \#\{i \mid t_i =1 \}$.  
If $s=1$, then we can take $V=U'$ and $J_t=L'_t$. Hence $\delta=0$ follows from Theorem \ref{split} directly. 
On the other hand, if $s>1$, then applying Theorem \ref{split} to $U$ and $U'$ we get the following two short exact sequences:
$$0\to H^1(V, J_t) \to H^1(U, L_t) \to  \oplus_{i=1}^s H^0(C_i^*, L''_t(i)) \to 0 $$
$$ 0\to H^1(V, J_t) \to H^1(U', L'_t) \to  \oplus_{i=2}^s H^0(C_i^*, L''_t(i)) \to 0$$
Putting these two short exact sequences together, we get that 
$$ 0\to H^1(U', L'_t) \to H^1(U, L_t) \to   H^0(C_1^*, L''_t(1)) \to 0$$
Hence the boundary map $\delta$ is always trivial.
\end{proof}

\br A similar proof is also used for \cite[Theorem 1.5]{BL} to study the length of the perverse sheaves $Rg_* L_t[2]$, where $g$ is the open inclusion form $U$ to $B$. The proof of Theorem \ref{main} indeed shows that the length of the perverse sheaves $Rg_* L_t[2]$ satisfies the formula in \cite[Theorem 2.4]{BL} as equality for any $t\in \cha(U)$, see \cite[Remark 2.5(a)]{BL}.  As mentioned in \cite[section 5, Remark 2.5(d)]{BL}, once such kind of formula holds as equality, one can get a formula for the characteristic cycles of the perverse sheaves $g_{!*} L_t[2]$. When $t_i=t_j$ for any $1\leq i,j \leq r$, this formula is not new, and it coincides with the one computed in \cite[Theorem 1.1]{NT}.
\er

%%%%%%%%%%%%%%%%%%%%%%%%%%%%%%%%%%%%%

%%% ====================== End of main part ====================== %%%
%\newpage

%------------------------------------------------------------------
\end{document}